# Properties of Commutative Association Schemes derived by FGLM Techniques.[*]

Edgar Martínez-Moro [†]
Dpto. Matemática Aplicada Fundamental
Universidad de Valladolid.
Valladolid, Castilla, 47002 Spain
`edgar.martinez@ieee.org`

January 24, 2000

**Abstract**

Association schemes are combinatorial objects that allow us solve problems in several branches of mathematics. They have been used in the study of permutation groups and graphs and also in the design of experiments, coding theory, partition designs etc. In this paper we show some techniques for computing properties of association schemes. The main framework arises from the fact that we can characterize completely the Bose-Mesner algebra in terms of a zero-dimensional ideal. A Gröbner basis of this ideal can be easily derived without the use of Buchberger algorithm in an efficient way. From this statement, some nice relations arise between the treatment of zero-dimensional ideals by reordering techniques (FGLM techniques) and some properties of the schemes such as P-polynomiality, and minimal generators of the algebra.

## 1 Overview and background.

In this paper the author is interested in the development of the connection between the algebraic properties of the Bose-Mesner algebra associated to an association scheme and the combinatorial properties of the last one. This investigation is mainly motivated by the fact that computing the eigenvalues of an association scheme is in general difficult (even numerically). We are interested in applications of the theory of association schemes to the areas of coding theory and design theory [4, 3, 7, 8, 15, 20], therefore, we mainly treat commutative symmetric association schemes. Anyway we do not discard the non commutative case and some directions are showed in the conclusions. Main tools the author proposes to show the above connection is using some reordering techniques from zero-dimensional ideals called FGLM techniques. The main objective is pointing the connection between the computer algebra treatment of zero-dimensional ideals and association schemes. In fact much about metric properties (P-polynomial properties) and Q-properties can be shown just by considering an appropriate set of equations. The paper is organized as follows: in the first section basic facts on association schemes as well as FGLM techniques are dealt with. We mainly follow the definition of an association scheme given in [1]. For an account on association schemes see for instance [15, 17, 25]. The second section establish the connection mentioned above and work out some relations. Some topics of these parts can be found in [21]. In the third section we build a correspondence between properties found by FGLM techniques and combinatorial properties of association schemes. Most of the new material is placed in this last part, even though some approaches to known combinatorial results in section 2 are also new. Finally in the conclusions section we show some open questions and lines of research. The approach in this paper is a computational one, as a final remark to this introduction we point that some other computational approaches can be found in [18] motivated by the works in [17, 19].

---

[*] To appear in *International Journal of Algebra and Computation*.
Preprint 88, Department of Algebra and Geometry (University of Valladolid)
[†] Partially supported by Dgicyt PB97-0471.

## 1.1 Association schemes.

Throughout this paper we will follow the definition of association scheme given by **E. Bannai and T. Ito** [1]:

**Definition 1.** An association scheme with $d$ classes is a pair $\mathcal{S} = (X, \{R_i\}_{i=0}^d)$, given by a finite set $X$ and a set of relations $\{R_i\}_{i=0}^d$ on $X$, satisfying the following rules:

1. $R_0 = \{(x, x) \mid x \in X\}$ (the diagonal relation)

2. $\{R_i\}_{i=0}^d$ is a partition on $X \times X$.

3. $\forall i \in \{0, \ldots d\} \quad \exists j \in \{0, \ldots d\}$ such that $R_i^t = R_j$, where
$$R_i^t = \{(y, x) \mid (x, y) \in R_i\}$$

4. For each election of $i, j, k \in \{0, \ldots d\}$, the number:
$$p_{ij}^k = |\{z \in X \mid (x, z) \in R_i \quad (z, y) \in R_j\}|$$
is constant for all $(x, y) \in R_k$

We can rewrite the above conditions as matrix relations. Consider the set of square matrices of order $v = |X|$ given by:

$$i = 1, \ldots, d \qquad D_i = [D_i(x, y)]_{0 \leq x, y \leq v} \qquad D_i(x, y) = \begin{cases} 1 & \text{if } (x, y) \in R_i \\ 0 & \text{elsewhere} \end{cases}$$

Therefore, the conditions $1 - 4$ above are equivalent to:

1'. $D_0 = Id$ (identity matrix)

2'. $\sum_{k=0}^d D_k = J$, where $J$ is the matrix where all entries are $1$.

3'. $\forall i \in \{0, \ldots d\} \quad \exists j \in \{0, \ldots d\}$ such that $D_i^t = D_j$

4'. $D_i \cdot D_j = \sum_{k=0}^d p_{ij}^k D_k$

We say that the scheme is commutative if $p_{ij}^k = p_{ji}^k \quad \forall i,j,k \in \{0,\ldots d\}$. We say the scheme is symmetric if $R_i^t = R_i$. The set of matrices $\{D_i\}_{i=0}^d$ is the generating set of a semisimple algebra $\mathcal{B}$ over $\mathbb{C}$ called **Bose-Mesner algebra** or **adjacency algebra**. All matrices in the set are linearly independent by condition $2'$, hence $\mathcal{B}$ has dimension $d+1$. If the scheme is commutative, it is clear that $D_i D_j = D_j D_i$, and therefore the algebra is commutative. From now on, for the sake of brevity, we will write scheme for a commutative symmetric association scheme.

Let $I$ be the identity matrix and $J$ the $v \times v$ matrix where all the entries are 1. It is well known (see [1]) that $\mathcal{B}$ is a diagonalizable algebra[1] and there is a unique set of primitive idempotents $\{E_0 = \frac{1}{|X|}J, E_1, \ldots, E_d\}$.

**Definition 2.** *The matrix representing the change of base from the base of the adjacency matrices to the idempotent matrices is called* **character table**[2] *and it is denoted by:*

$$P = [p_i(j)]_{i,j=0}^d \tag{1}$$

Note that the entries on the $i^{th}$ row are the eigenvalues (counting multiplicities) of the matrix $D_i$. Also, since we are considering symmetric association schemes, the $D_i$'s are symmetric and the eigenvalues are real. If we consider the algebra $\mathcal{U}$ generated by $\{E_0, E_1, \ldots, E_d\}$ endowed with Hadamard product over the complex numbers, we have that the set of primitive idempotents is the set $\{D_i\}_{i=0}^d$ (see [1]) and:

$$E_j = \sum_{i=0}^d q_i(j) D_i \qquad Q = [q_i(j)]_{i,j=0}^d = |X|P^{-1} \tag{2}$$

The $q_i(j)$'s are called *second eigenvalues*.

## 1.2 Gröbner bases techniques.

In this paper we use some well-known facts on Gröbner bases theory, for background material on this theory see [6] or [24]. Main topics we need is the relation between diagonalization and solving 0-dimensional ideals and some clues on FGLM results. Treatment of finite generated algebras can be found in [9, 14, 23]

### 1.2.1 Solving 0-dimensional ideals and Gröbner bases.

First we show some relations of the problem of solving a 0-dimensional system of equations with an eigenvalue problem. This material can be found in [22, 5]. Indeed this fact discloses the relation between the systems proposed in the next section and the problem of diagonalizing the Bose-Mesner algebra. Roughly speaking it works as follows:

We are given an ideal $\mathcal{I}$ generated by $\mathcal{F}$, where $\mathcal{F}$ is a reduced Gröbner basis for $\mathcal{I}$ in some monomial ordering and $V(\mathcal{F})$ (the variety generated by $\mathcal{F}$) is zero-dimensional. As usual $lt(f)$ will denote the leading term of $f$ and the set of terms will be denoted by:

$$T = \{x_0^{i_0} \ldots x_d^{i_d} \mid i_0, \ldots, i_d \in \mathbb{N}_0\} \tag{3}$$

---

[1] This arises from the fact that $\{D_i\}$ is a set of commuting matrices.
[2] The justification for this name is given in [25].

We define the normal set as:

$$N = \{t \in T \mid \nexists f \in \mathcal{F} \text{ such that } lt(f)|t\} \tag{4}$$

$N$ is a basis of $\mathbb{C}[x_0, \ldots, x_d] \bmod \mathcal{I}$ and the cardinality of $N$ is just the number of roots of $\mathcal{F}$. We can describe the effect of multiplying an arbitrary element by a fixed $f \in \mathbb{C}[x_0, \ldots, x_d]$ modulo $\mathcal{I}$ just by using the images of multiplying $f$ by each term in $N$. Hence, if we denote by $\mathrm{nf}(f)$ the normal form with respect $\mathcal{F}$, we conclude that:

$$\mathrm{nf}(f \cdot x_i) = \sum a_{ij}(f) x_i \quad i = 0, \ldots d \tag{5}$$

As usual the complex matrix $A(f) = [a_{ij}(f)]_{i,j=0}^d$ is the multiplication matrix of $f$. The multiplication tables for the elements in $N$ are particularly important since from them we can rebuild any multiplication table and clearly $A(x_i) = [p_{ij}^k]_{j,k=0}^d$. The following result is a slightly modified version of a theorem in [5, 22].

**Theorem 1 (Möller, Stetter).** *Let $\mathcal{F}$ a Gröbner basis, $M_0, \ldots, M_d$ the multiplication matrices for $x_0, \ldots, x_d$ defined as above. Let $U$ the unitary matrix such that $U^* M_i U = F_i$ is diagonal for each $i \in \{0, 1, \ldots, d\}$. If we denote the diagonal entries of the matrices $F_i$ as $(u_0^i, \ldots, u_d^i)$ then the points on $V(\mathcal{F})$ counting multiplicity are:*

$$z_j = (x_0^j, \ldots, x_d^j) = (u_j^0, u_j^1, \ldots, u_j^d) \tag{6}$$

*i.e. the roots of $\mathcal{F}$ are the eigenvalues of the $M_i$.*

*Proof.* For a proof of the result see [22]. □

Therefore, as usual, we will reduce the problem of diagonalizing the algebra $\mathcal{B}$ to solve a system of equations in $d$ indeterminates over $\mathbb{C}$.

### 1.2.2 FGLM results.

Finally in this section we point here briefly to two properties on Gröbner bases for 0-dimensional ideals. An introductory text to these ideas can be seen in [6] or [24], the techniques and results exposed here can be found in [13, 12]. We say that an ideal $I \in \mathbb{R}[x_0, \ldots, x_n]$ is in generic position with respect $n$ if, for each $(a_0, \ldots, a_n)$ and $(b_0, \ldots, b_n)$ distinct zeros of $I$ in $\mathbb{R}^n$, $a_n \neq b_n$. We say that a change of coordinates is a generic change of coordinates for an ideal $I$ if it is in generic position for one of the variables.

**Proposition 1 (Gianni,Mora).** *Let $I$ be a zero-dim. ideal in $\mathbb{R}[x_0, \ldots, x_n]$, there is a generic change of coordinates for $I$, and there are $g_0, \ldots, g_n \in \mathbb{R}[x_n]$ such that if $i < n$ then $degree(g_i) < degree(g_n) = d$ where $d$ is the number of distinct zeros of $I$ and $g_n$ is squarefree, and finally*

$$\mathcal{G} = \{x_0 - g_0, x_1 - g_1, \ldots, x_{n-1} - g_{n-1}, g_n\}$$

*a reduced Gröbner basis of $I$ w.r.t. the lexicographical ordering such that $x_0 > x_1 > \cdots > x_n$*

*Proof.* See the proof of the first assumption in proposition 1.6 in [12]. □

**Proposition 2.** *There is a probabilistic algorithmic method of computing such a generic change of coordinates.*

*Proof.* For a proof see [12], in fact, that paper is devoted to compute such an expresion. □

We will show the combinatorial significance of these facts in the case of the Bose-Mesner algebra later in this paper.

## 2 Ideals related to association schemes.

In this section we introduce some systems of equations having as solutions the first and second eigenvalues of the scheme. We will interpret some results both in terms of computer algebra and combinatorics, showing the great resemblance between both topics.

### 2.1 First eigenvalues. Interpretation of the intersection algebra.

It is a well known fact in combinatorics that the first eigenvalues of the scheme defined above are the nonzero solutions of the system given by:

$$x_i x_j - \sum_{k=0}^{d} p_{ij}^k x_k = 0 \quad 0 \leq i \leq j \leq d \tag{7}$$

This fact can be proved in many ways, we will sketch a proof in a Gröbner basis flavor. For a complete and rigorous proof see [21], here we will be more concerned with showing the parallel noted above.

Let us consider the set of polynomials (as usual $x_0 = 1$, i.e. the identity relation):

$$\mathcal{F} = \{x_i x_j - p_{ij}^0 - \sum_{k=1}^{d} p_{ij}^k x_k\}_{1 \leq i \leq j \leq d} \tag{8}$$

From now on we will call equations in (7) structure equations and the polynomials in $\mathcal{F}$ structure polynomials. There are some nice properties on $\mathcal{I}$ the ideal generated by $\mathcal{F}$:

1. The variety $V(\mathcal{I})$ is just the set of nonzero solutions of the equations in (7).

2. $\mathcal{F}$ is a reduced Gröbner basis for $\langle \mathcal{F} \rangle = \mathcal{I}$ in the total degree monomial ordering.

3. $\mathcal{I}$ is a radical ideal, ie. the variety $V(\mathcal{I})$ is reduced.

We recall the above section and build up the multiplication matrices $A(f) = [a_{ij}(f)]_{i,j=0}^{d}$ for $f \in \mathcal{F}$. The multiplication tables for the elements in the corresponding normal set $N = \{1, x_1, \ldots, x_d\}$ are the matrices given by $M^i, 0 \leq i \leq d$ where:

$$M^i = \left(m_{lm}^i\right)_{l,m=0}^{d}, \quad \text{where } m_{lm}^i = p_{li}^m \tag{9}$$

These matrices are known in combinatorics as **intersection matrices** and it is well known that they generate an algebra isomorphic to the Bose-Mesner algebra called intersection algebra. Solving the equations in (7) is the same as diagonalising the intersection algebra and, since they are isomorphic, is the same as diagonalising the Bose Mesner algebra. This is the main fact that offers us a bridge between properties of these algebras and properties of the ideal $\mathcal{I}$.

## 2.2 On the other hand ...

Another well known fact is that the solutions of the $d$ quadratic equations:

$$x_k - \sum_{j=1}^{d} p_{0j}^k x_k - \sum_{i,j=1}^{d} p_{ij}^k x_i x_j = 0 \quad 0 \leq k \leq d \tag{10}$$

is the variety of idempotents or projectors of the algebra. I.e. the solutions are the second eigenvalues of the scheme. The main trade off of this system is that the polynomials involved do not form a Gröbner basis of the variety and in general, computing Gröbner basis is a non-efficient problem, but once achieved a Gröbner basis of it, all computations are carried out in the same way. Some tricks to deal with these equations are sumarized in the reference [24]. We will consider only the first eigenvalues from now on.

# 3 Combinatorial properties achieved by reordering techniques.

One of the main facts when dealing with association schemes is to know whether they are P-polynomial or not. We say that a scheme $\mathcal{S} = (X, \{R_i\}_{i=0}^{d})$ is P-polynomial if we can reorder the relations $R_i$ so that the corresponding $D_i$ is a polynomial $p_i$ in $D_1$ with degree $i$ and $p(x_1)$ is a polynomial in $x_1$ whose roots are the eigenvalues of $D_1$.

Clearly, if $\mathcal{S}$ is a P-polynomial scheme we have that there is a monomial ordering such that the set:

$$\{x_0 - 1, p(x_{i_1}), x_{i_2} - p_2(x_{i_1}), \ldots x_{i_d} - p_d(x_{i_1})\} \tag{11}$$

is a reduced Gröbner basis for $\mathcal{I}$ for a pure lexicographical order where $x_{i_1} < x_{i_j}$, $j = 2, \ldots, d$. It follows that next proposition holds:

**Proposition 3.** *Let $\mathcal{S} = (X, \{R_i\}_{i=0}^{d})$ be an association scheme and $\mathcal{I}$ its associated ideal. $\mathcal{S}$ is a P-polynomial scheme if and only if there is an pure lexicographical ordering for the variables $x_{i_j}$ with $x_{i_0} = x_0, x_{i_1} < x_{i_j}, j \neq i_0$ such that*

$$\{x_0 - 1, p(x_{i_1}), x_{i_2} - p_2(x_{i_1}), \ldots x_{i_d} - p_d(x_{i_1})\} \tag{12}$$

*is a Gröbner basis for $\mathcal{I}$ for that ordering, and the following conditions on the degree of the polynomials $p, p_1, \ldots, p_d$ are satisfied:*

$$degree(p_i) = i \quad \forall i = 2, \ldots, d \quad degree(p) = d \tag{13}$$

Note that for an association scheme to be P-polynomial its Bose-Mesner algebra is the minimal algebra containing one of the relations but the converse is not true in general, see example 19 in [21]. The reason for this is the conditions on the degrees of the polynomials. Note also that the resulting univariate polynomial is squarefree and it is just the characteristic polynomial of the relation that "spans" the whole algebra[3]. The main profit we get from proposition above is an algorithmic way of checking the P-polynomiality of an association scheme:

**Algorithm 1 (Checking P-polynomiality).**

- **Input** $\mathcal{F}$ Structure equations of $\mathcal{S}$,
- **For** $i$ from $1$ to $d$ **do**
  **Compute** a reduced Gröbner basis for $\langle \mathcal{F} \rangle$ for a pure lexicographical monomial ordering where $x_i < x_j \; j \neq i$
  **Check** the conditions on proposition above.
  **If** the conditions are fulfilled then **Stop** and
  **Return** The scheme is P-polynomial and the Gröbner basis computed **fi**
  **od**
  **Stop Return** The scheme is not P-polynomial

Complexity of the algorithm is placed in the Gröbner basis computation, but since we are dealing 0-dimensional ideals and $\mathcal{F}$ is already a Gröbner basis, these computations can be done by FGLM techniques which are effective and efficient. Clearly the complexity is at most the same as $d$ FGLM Gröbner basis computations. The algorithm can be easily modified to check whether the scheme is P-polynomial with respect to more than one relation or to check the minimal number of relations spanning the whole algebra. This last fact is important in coding theory for computing the Lloyd polynomials and metric properties of the scheme [21]. In that paper can be found many examples of how to apply this algorithm, we recall just two of them:

**Example 1.** Let $X = \mathbb{Z}_m, m \in \mathbb{N}$ and a radix $r$ $(r < m)$. The subgroup of $\mathbb{Z}_m^*$ generated by $r$ and $-1$ acts on $X$ by multiplication with orbits $\mathcal{O}_k$. The following relation defines an association scheme (2-orbit association scheme) on $X$:

$$xR_k y \Leftrightarrow x - y \in \mathcal{O}_k$$

When $m = 9, r = 2$ the orbits are: $\mathcal{O}_0 = \{0\}, \mathcal{O}_1 = \{-1, 1 - 2, 2, -4, 4\} \, \mathcal{O}_2 = \{-3, 3\}$ $\mathcal{F}$ is given by:

$$\mathcal{F} = \{x_2{}^2 - 6 - 6x_1 - 3x_2, \; x_0 - 1, \; x_1{}^2 - x_1 - 2, \; x_1 x_2 - 2x_2\}$$

The above algorithm return us the following Gröbner basis:

$$\{-x_2{}^2 + 6 + 6\,x_1 + 3\,x_2, \; x_0 - 1, \; -18\,x_2 - 3\,x_2{}^2 + x_2{}^3\}$$

Which allows us to compute the character table by solving the equation $-18\,x_2 - 3\,x_2{}^2 + x_2{}^3 = 0$ and evaluating the solutions in the other expressions:

$$CharTab = \begin{pmatrix} 1 & 2 & 6 \\ 1 & 2 & -3 \\ 1 & -1 & 0 \end{pmatrix}$$

---

[3]This means that any element in the algebra can be expressed by a formula of the generator and $J$, using the addition, the ordinary multiplication and the Hadamard product.

The following example is not P-polynomial:

**Example 2.** *Consider the scheme defined as above given by $m = 8, r = 3$, the orbits are: $\mathcal{O}_0 = \{0\}, \mathcal{O}_1 = \{\pm 1, \pm 3\}, \mathcal{O}_2 = \{\pm 2, \}, \mathcal{O}_3 = \{4\}$ The associated ideal $\mathcal{I}$ is generated by:*

$$\{x_1{}^2 - 4x_2 - 4x_3 - 4, x_3{}^2 - 1, x_0 - 1, x_1 x_2 - 2x_1,$$
$$x_1 x_3 - x_1, x_2{}^2 - 2x_3 - 2, x_2 x_3 - x_2\}$$

*If we compute the reduced Gröbner basis with respect plex. orders we have:*

1. $x_1 < x_2, x_3$

   $4x_3 + 4x_2 - x_1{}^2 + 4, 2x_2{}^2 + 4x_2 - x_1{}^2, x_1 x_2 - 2x_1, -16x_1 + x_1{}^3, x_0 - 1$

2. $x_2 < x_1, x_3$

   $x_0 - 1, x_1{}^2 - 4x_2 - 2x_2{}^2, x_1 x_2 - 2x_1, -x_2{}^2 + 2x_3 + 2, -4x_2 + x_2{}^3$

3. $x_3 < x_1, x_2$

   $x_1{}^2 - 4x_2 - 4x_3 - 4, x_1 x_2 - 2x_1, x_1 x_3 - x_1, x_2{}^2 - 2x_3 - 2, x_2 x_3 - x_2,$
   $x_3{}^2 - 1, x_0 - 1$

*None of the three cases above fulfill the conditions on the proposition, and therefore is not P-polynomial.*

In fact, as we said above, a small modification on the algorithm allows us to compute the minimal number of relations spanning the algebra:

**Example 3.** *Following the setting in the example above, if we choose one ordering (for example $x_1 < x_2, x_3$ ) we have the following expression in two variables of the relations:*

$$x_0 = 1, x_1 = x_1, x_2 = x_2, x_3 = \frac{1}{4}\left(x_1{}^2 - 4x_2 + 4\right)$$

Finally, propositions 1 and 2 allows us find an algorithmic procedure to find a non-negative integer matrix $A$ in the Bose Mesner algebra such that any element in the algebra is *expressible in terms of* $A$, i.e. spans the whole algebra[4]. Clearly, this happens if and only if there is a change of coordinates such that $A$ is a resulting coordinate, and $I$ is in generic position w.r.t. $A$. Therefore lemma 1.1.1 in [17] becomes obvious in a Gröbner basis setting also as they recall in their paper from a linear algebra approach. This matrix is of great importance when we compute partition designs on the association scheme. A modification of algorithm 5.4 in [12] allows us find such an element $A$, the idea is making changes of coordinates until the ideal is in generic position, this is achieved by proposition 2. For the understanding of the modifications we made in the algorithm one must take in mind that our ideal is already radical, so we do not need to compute the radical, and also that we are looking for an non-negative integer change of coordinates which imposes no restriction on the running of the algorithm (For a proof see [13]).

---
[4]The meaning of this matrix $A$ in terms of graph theory can be found in [17]

**Algorithm 2 (Finding the matrix $A$).**

- **Input** $\mathcal{F}$ Structure equations of $\mathcal{S}$,
- **Compute** a reduced Gröbner basis $\mathcal{G}$ for $\langle \mathcal{F} \rangle$ for a pure lexicographical monomial ordering where $x_d < x_j$ $0 \leq j < d$
- **Repeat**
  **Let** $\{g\} = \mathcal{G} \cap \mathbb{R}[x_d]$.
  **Let** $i$ the maximal index s.t. for $j \leq i$, $x_j$ is expressible in terms of $x_d$.
  **If** $i < d-1$ **then**
  **Choose a random** positive integer $c$
  **Make** the change of coordinates[5]

$$x_j \mapsto x_j \quad 0 \leq j < d$$
$$cx_i + x_d \mapsto x_d$$

  **Compute** a reduced Gröbner basis $\mathcal{G}$ with the same ordering in the new variables. **fi**
  **Until** $i = d-1$ **Stop**
- **Return** The Gröbner basis and the change of coordinates to achieve the last $x_d$.

Readers interested in the algorithm are referred to the paper [12] for a more detailed description of why it works. We just show how to deal with example 2 above:

**Example 4.** *Taking the change of coordinates given by:*

$$x_i \mapsto y_i \quad 0 \leq i < 3$$
$$x_3 + x_1 \mapsto y_3$$

*The structure equations in the new base become:*

$$y_1^2 - 4y_2 - 4(y_3 - y_1) - 4,\ (y_3 - y_1)^2 - 1,\ y_0 - 1,\ y_1 y_2 - 2y_1$$
$$y_1(y_3 - y_1) - y_1,\ y_2^2 - 2(y_3 - y_1) - 2,\ y_2(y_3 - y_1) - y_2$$

*And a Gröbner basis for a plex. ordering where $y_3 < y_1, y_2$ is given by:*

$$12y_1 - y3^3 + 3y_3^2 + y_3 - 3,\ 27 + 24y_2 - 3y_3^2 + 25y_3 - y_3^3,$$
$$15 - 16y_3^2 + 2y_3 - 2y_3^3 + y_3^4,\ x_0 - 1$$

*Hence the ideal is in generic position and the algebra is expressible in terms of the matrix given by the sum of the relation $D_1$ and relation $D_3$.*

## 4 Conclusions

Properties of symmetric commutative association schemes can be derived from their associated algebras via ideal computations and Gröbner basis of the ideals $\mathcal{I}$ associated with then allow us to describe a dictionary between combinatorial properties and algebraic treatment of zero-dimensional ideals (see appendix). We have proved that these ideals are radical, and therefore any polynomial condition on the points $V(\mathcal{I})$ can be seen as a membership problem which is answered by computing the normal form. All computations and algorithms proposed are effective and efficient since we are dealing with zero-dimensional Gröbner basis using FGLM techniques [10] on the ideal. It is not surprise one can use Gröbner methods here since we are dealing with finite dimensional algebras, and the linear algebra in the background is indeed that those researchers working on association schemes use. The main contribution of this paper is show this connection and also the algorithms on the paper give some extra information about the association scheme for free. Algoritm 1 does not only check P-polynomiality but gives us the minimal number of relations generating the algebra and the expresion of the Lloyd's polynomials and the matrix $A$ resulting from algoritm 2 gives us a valuable tool for checking partition designs on the scheme.

Moreover, future trends of work point towards extending this setting to the non-commutative case using the FGLM techniques in [2]. We also would like to study the relationship pointed by the Galois correspondence studied in [18] and the corresponding ideals. Finally a new direction of work is develop a similar approach for dealing with proper and adjacency polynomials in the theory of graph spectra [11].

# 5  Appendix.

This paper allow us compose a small dictionary translating combinatorial properties of the association scheme and properties and algorithms applied to is associated 0-dimensional ideal. It can be set down as follows:

| Association Schemes | Algebraic theory of finite algebras. |
|---|---|
| $\mathcal{S} = (X, \{R_i\}_{i=0}^d)$ | $\mathbb{C}[x_0, \ldots, x_d]/\langle \mathcal{F} \rangle$ |
| Character table of $\mathcal{S}$ | Points on $V(\langle \mathcal{F} \rangle)$ |
| P-Polynomiality | Algorithm 1. |
| Minimal $\sharp$ of relations expressing the scheme | Minimal $\sharp$ of generator for the algebra Algorithm 1 modified. |
| Single matrix for expressing the scheme | General change of coordinates for $\langle \mathcal{F} \rangle$ Algorithm 2. |
| Metric properties | [21] |

# 6  Acknowledgements.

The author would like to thank Dr. A. Campillo-López for his reading of the preliminary manuscript and to Professor Michail Klin for the availability of all his material about this topic. All examples (orbits, Gröbner basis, etc.) have been computed using Maple V release 4.

# References

[1] **E. Bannai, T. Ito** Algebraic combinatorics I: association schemes. *Benjamin-Cumming publishers*, 1984.

[2] **M.A. Borges, M. Borges, T. Mora** Computing Gröbner bases by FGLM Techniques in a Noncommutative Setting, *J. Symb. Comp*, To appear.

[3] **P.J. Cameron, J.H. Van Lint** Graphs, codes and designs. *Cambridge University Press* (LMS lecture notes 43) 1980.